\documentclass[12pt,a4paper]{amsart}
\usepackage{amssymb}

\usepackage{amsmath}
\usepackage{amsfonts}

\newcommand{\Z}{\mathbb{Z}}
\newcommand{\F}{\mathbb{F}}

\DeclareMathOperator{\Aut}{Aut}

\newtheorem{main-dummy}{Main-Dummy}

\newtheorem{main-theorem}[main-dummy]{Theorem}

\numberwithin{dummy}{section}
\numberwithin{equation}{section}

\newtheorem*{lemma*}{Lemma}

\newtheorem*{theorem*}{Theorem}

\newtheorem*{prop*}{Proposition}

\theoremstyle{definition}

\newtheorem*{example*}{Example}

\theoremstyle{remark}

\newtheorem*{rem*}{Remark}

\begin{document}
\bibliographystyle{amsalpha}
\author{Sandro Mattarei}

\email{mattarei@science.unitn.it}

\urladdr{http://www-math.science.unitn.it/\~{ }mattarei/}

\address{Dipartimento di Matematica\\
  Universit\`a degli Studi di Trento\\
  via Sommarive 14\\
  I-38050 Povo (Trento)\\
  Italy}

\title{A note on automorphisms of free nilpotent groups}

\begin{abstract}
We exhibit normal subgroups of a free nilpotent group $F$ of rank two and class three,
which have isomorphic finite quotients but are not conjugate under any automorphism of $F$.
\end{abstract}

\subjclass[2000]{Primary 20E05; secondary 20F18}

\keywords{Free nilpotent group, automorphism}

\thanks{The  author  is grateful  to  Ministero dell'Istruzione, dell'Universit\`a  e
  della  Ricerca, Italy,  for  financial  support of the
  project ``Graded Lie algebras  and pro-$p$-groups of finite width''.}

\maketitle

\thispagestyle{empty}

A remarkable fact about free profinite groups of finite rank
is that any isomorphism between finite quotients of such a group $F$ lifts to an automorphism of $F$.
This is true, more generally, if $F$ a free pro-$\mathcal{C}$-group of finite rank,
where $\mathcal{C}$
is a family of finite groups closed under taking subgroups, homomorphic images and direct products,
and containing nontrivial groups.
A proof and the relevant definitions can be found in~\cite[Proposition~15.31]{FriedJarden},
but the essence of the argument goes back to Gasch\"utz~\cite{Gaschutz}.
In preparation for a summer school on {\em ``Zeta functions of groups''}
held by Marcus du Sautoy and the author in June 2002 in Trento (Italy),
du Sautoy suggested that this may remain true
for (abstract) free nilpotent groups $F$, and asked the author,
who was responsible for that part of the course, to provide a proof.
If confirmed, this claim would have simplified the course
by avoiding the need to set up the language of profinite groups.

Unfortunately, this claim already fails for $F$ a free abelian group of rank one, that is,
an infinite cyclic group:
in this case $F$ has exactly two automorphism, while its quotient of order $n$
has $\varphi(n)$ automorphisms, and $\varphi(n)>2$ for $n>4$.
A milder statement which would have been sufficient for our purposes
would be that any two normal subgroups of $F$ with finite isomorphic quotients
are conjugate under some automorphism of $F$.
This is also false, and one does not have to dig much deeper in order to find a counterexample.
We first record an example suggested by the anonymous referee.
It is based on a three-generated group of order $p^6$ and class two,
which was studied in~\cite{DauHei}.
After that we present an example where $F$ is two-generated
and the quotients have order $p^4$.
This order is easily seen to be minimal for such an example.

\begin{example*}
The groups $G$ of odd order $p^6$ satisfying
$G'=Z(G)=G^p$ were classified by Daues and Heineken in~\cite{DauHei}
in terms of dualities of a three-dimensional vector space over the field of $p$ elements.
In particular, the group $G$ in their case (I) has a $p$-group as the full group of automorphisms.
One can realize $G$ as the quotient of the free nilpotent group $F=\langle x,y,z\rangle$
of rank three and class two modulo the normal subgroup
$M_r$ generated by $(F')^p$ and the three elements
\[
x^{rp}[y,x],\quad
y^{rp}[z,x],\quad
z^{rp}[z,x]^{-1}[z,y],
\]
where $r$ is any integer prime to $p$.
When $r=1$ the relations associated with these three elements correspond to the matrix
$D=(a_{ij})=
\left(\begin{smallmatrix}
0&0&1\\0&-1&0\\1&1&0
\end{smallmatrix}\right)$
as described in~\cite[p.~219]{DauHei}, with $x,y,z$ in place of $x_1,x_2,x_3$.
However, all choices of $r$ prime to $p$ give rise to isomorphic groups $F/M_r$.
Assuming $p\neq 2,3,7$, we can choose $r$ such that $r^3\not\equiv\pm 1\pmod{p}$.
In particular, we may always take $r=2$.
Then $M_1$ and $M_r$ are not conjugate under $\Aut(F)$.

In order to see this it suffices to show that no isomorphism
of $F/M_r$ onto $F/M_1$ lifts to an automorphism of $F$.
One isomorphism of $F/M_r$ onto $F/M_1$ is obtained by mapping
$x,y,z$ to $x^r,y^r,z^r$, respectively.
This induces an automorphism of their common quotient $F/M_1M_r=F/F'F^p$,
with determinant $r^3$ when the latter is viewed as a vector space over the field of $p$ elements.
Any other isomorphism of $F/M_r$ onto $F/M_1$ is obtained by composing the one described
with an automorphism of $F/M_1$.
Since the latter has $p$-power order, and hence determinant one on $F/F'F^p$,
we conclude that every isomorphism of $F/M_r$ onto $F/M_1$ induces an automorphism of $F/F'F^p$ with determinant $r^3$.
Because automorphisms of $F$ induce maps of determinant $\pm 1$ on
$F/F'$ viewed as a free $\Z$-module, and $r^3\not\equiv\pm 1\pmod{p}$, they cannot induce any
isomorphism of $F/M_r$ onto $F/M_1$.

A more careful analysis, such as that in the proof of the Theorem below, would reveal that
for any odd prime $p$ (thus including $3$ and $7$), the subgroups
$M_r$ and $M_s$ are conjugate under $\Aut(F)$ if and only if $r\equiv\pm s\pmod{p}$.
We leave the details to the interested reader and only suggest to use the fact that
the subgroups
$\langle G',x\rangle$ and
$\langle G',x,y\rangle$
of $G=F/M_1$ are characteristic.
In fact, according to~\cite{DauHei},
$\Aut(G)$ is generated by
the automorphism determined by
$x\mapsto x$,
$y\mapsto xy$,
$z\mapsto yz$
together with the $p^9$ central automorphisms, which induce the identity map on $G/G'$.
\end{example*}

In the two-generated example which we present now the group of automorphisms of the finite quotients
is not a $p$-group.
Hence the proof is more involved, and we formally state the result as a theorem.

\begin{theorem*}
Let $F=\langle x,y\rangle$ be the free nilpotent group of rank two and class three,
and let $p$ be a prime greater than three.
For $r=1,\ldots,p-1$ set
\[
N_r=\langle x^{p^2},y^p,x^{-rp}[y,x,x],[y,x,y]\rangle^F,
\]
where the superscript $F$ denotes taking the normal closure in $F$.
Then $F/N_r$ is a $p$-group of order $p^4$, class three and exponent $p^2$.
All quotients $F/N_r$ are isomorphic.
However, $N_r$ and $N_s$ belong to the same orbit under $\Aut(F)$ if and only if
$r=s$ or $r=p-s$.
\end{theorem*}

\begin{proof}
It is well known
that each element of $F$ can be written as
$x^iy^j[y,x]^k[y,x,x]^l[y,x,y]^m$
for uniquely determined integers $i,j,k,l,m$,
see~\cite[Theorem~11.2.4]{MHall}.
It is then easy to deduce that each coset of
$K=\langle x^{p^2},y^p,[y,x,y]\rangle^F$ in $F$ has a unique representative of the form
$x^iy^j[y,x]^k[y,x,x]^l$, with $0\le i<p^2$ and $0\le j,k,l<p$.
This also follows from a general result giving $\F_p$-bases, in terms of basic commutators
and their powers, for the factors of the lower $p$-central series of a free group,
see~\cite[Lemmas 1.11 and 1.12]{Scoppola:Frobenius-Wielandt}, for instance.
In particular, $K$ has index $p^5$ in $F$, and hence $N_r=\langle K,x^{-rp}[y,x,x]\rangle$
has index $p^4$ in $F$.
Clearly, $F/N_r$ has class three and exponent $p^2$.

We will determine all endomorphisms of $F$ which map
$N_r$ into $N_s$ and induce an isomorphism between the quotient groups $F/N_r$ and $F/N_s$.
Since $M/N_r$, where $M=\langle x^p,y\rangle^F$, is the only abelian maximal subgroup of $F/N_r$,
we may restrict our attention to endomorphisms which map $M$ into itself.
Thus, let $\psi$ be an endomorphism of $F$ such that
$\psi(x)= x^iy^jc$ and $\psi(y)=y^kd$, for some integers $i,j,k$ and some $c,d\in F'F^p$.
We may also assume that $i$ and $k$ are prime to $p$, because this is a necessary condition for inducing an isomorphism
of $F/N_r$ onto $F/N_s$ and, in particular, an automorphism of $F/F'F^p$.

As a special case of~\cite[Hilfssatz~III.10.9(b)]{Hup} or~\cite[Corollary~1.1.7(i)]{L-GMcKay},
if $a,b$ are elements of a $p$-group $G$ of class less than $p$, and if the normal closure of $b$
is abelian of exponent $p$, then $(ab)^p=a^p$.
Since
$\langle y,[y,x],[y,x,x],K\rangle/K$,
the normal closure of $yK$ in $F/K$, is abelian of exponent $p$,
and because of standard commutator identities, we have
\begin{align*}
&\psi(x^{p^2})=((x^iy^jc)^p)^p\equiv (x^{ip})^p\equiv 1\pmod{K}
\\
&\psi(y^p)=(y^kd)^p\equiv 1\pmod{K}
\\
&\psi([y,x,y])=[y^kd,x^iy^jc,y^kd]\equiv [y,x,y]^{ik^2}\equiv 1\pmod{K}.
\end{align*}
Thus, $\psi$ maps $K$ into itself.
Because of our assumption that $i$ and $k$ are prime to $p$,
it induces an automorphism of $F/F'F^p$, and hence an automorphism of $F/K$,
since the former is the Frattini quotient of the latter.
Finally, we have
\begin{align*}
\psi(x^{-rp}[y,x,x])&=(x^iy^jc)^{-rp}[y^kd,x^iy^jc,x^iy^jc]
\\
&\equiv x^{-irp}[y,x,x]^{i^2k}
\pmod{K}.
\end{align*}
Consequently, $\psi$ maps $N_r$ into $N_s$ if and only if $x^{-irp}[y,x,x]^{i^2k}$
equals a power of $x^{-sp}[y,x,x]$, that is,
if and only if $iks\equiv r\pmod{p}$.
If this condition is met, and it certainly can by a suitable choice of $i$ and $k$, then
$\psi$ induces an isomorphism of $F/N_r$ onto $F/N_s$, as desired.

It remains to see when the endomorphism $\psi$ of $F$ is an automorphism.
Recall that $F$, being a finitely generated nilpotent group, is hopfian,
that is, each surjective endomorphism of $F$ is an automorphism~\cite[Theorem~5.5]{MKS}.
Thus, $\psi$ is an automorphism if and only if it is surjective, that is, if and only if it induces
an automorphism of its Frattini quotient $F/\Phi(F)=F/F'$.
It follows that $\psi$ is an automorphism of $F$ if and only if $ik=\pm 1$.
Consequently, $N_r$ and $N_s$ belong to the same orbit under $\Aut(F)$ if and only if $r\equiv\pm s\pmod{p}$.
\end{proof}

\bibliography{References}

\end{document}